\documentclass[12pt]{amsart}

\usepackage{amssymb}
\usepackage[all]{xypic}

\newcounter{alphalistcntr}

\theoremstyle{plain}
\newtheorem{theorem}{Theorem}[section]

\newtheorem{corollary}[theorem]{Corollary}
\newtheorem{proposition}[theorem]{Proposition}

\theoremstyle{remark}

\newtheorem{remark}[theorem]{Remark}

\newtheorem*{note*}{Note}
\newtheorem*{remark*}{Remark}
\newtheorem*{example*}{Example}

\theoremstyle{definition}
\newtheorem*{definition*}{Definition}

\providecommand{\Gal}{\mathop{\rm Gal}\nolimits}
\providecommand{\Cl}{\mathop{\rm Cl}\nolimits}
\providecommand{\D}{\mathop{\rm D}\nolimits}
\providecommand{\K}{\mathop{\rm K}\nolimits}
\providecommand{\mod}{\mathop{\rm mod}\nolimits}
\providecommand{\Mat}{\mathop{\rm Mat}\nolimits}
\providecommand{\Mod}{\mathop{\rm Mod}\nolimits}

\title[Capitulation for locally free class groups]{Capitulation for locally free class groups of orders of group algebras over number fields}

\author{Cornelius Greither}
\address{Cornelius Greither\\
Fakult\"at f\"ur Informatik\\
Institut f\"ur theoretische Informatik und Mathematik\\
Universit\"at der Bundeswehr M\"unchen\\
85577 Neubiberg\\
Germany}
\email{cornelius.greither@unibw.de}

\author{Henri Johnston}
\address{Henri Johnston\\ 
St. Hugh's College\\
St. Margaret's Road\\
Oxford OX2 6LE\\
United Kingdom}
\email{henri@maths.ox.ac.uk}
\urladdr{http://www.maths.ox.ac.uk/$\sim$henri}

\thanks{Johnston was supported by a grant from the
 Deutscher Akademischer Austausch Dienst.}
\subjclass[2000]{Primary 11R65; Secondary 11R33}
\keywords{Locally Free Class Groups, Additive Galois Module Structure}
\date{9th July 2009}

\begin{document}
\maketitle

\begin{abstract}
We prove a capitulation result for locally free class groups of orders of group algebras
over number fields. This result allows some control over ramification and so as a corollary
we obtain an ``arithmetically disjoint capitulation result'' for the Galois module structure of rings of integers.
\end{abstract}

\section{Introduction}\label{intro}

Let $K/F$ be a finite Galois extension of number fields with Galois group $G$. 
A natural problem that arises in Galois module theory is that of determining whether
the ring of integers $\mathcal{O}_{K}$ is free as a module over the associated order
\[
\mathcal{A}_{K/F} := \{ x \in F[G] \, \mid \, x(\mathcal{O}_{K}) \subseteq \mathcal{O}_{K}\}.
\]
One tool that has been used with some success is that of the locally free class group 
(see \cite{frohlich-book} for an introduction).

Now let $F$ be a number field and $H$ be its Hilbert class field. It is well-known that
every ideal of $\mathcal{O}_{F}$ becomes principal or ``capitulates'' in $\mathcal{O}_{H}$. 
In other words, the natural map $\Cl(\mathcal{O}_{F}) \rightarrow \Cl(\mathcal{O}_{H})$ is trivial. 
The main theorem of this paper is an analogous result for certain locally free class groups 
that allows some control over ramification 
(though unlike the case of the Hilbert class field, our construction is far from canonical). 
As a corollary, we obtain an ``arithmetically disjoint capitulation result'' for the Galois module 
structure of rings of integers (see \cite[III.2.13]{ft} for the definition of arithmetical disjointness).
Our results are precisely stated as follows.

\begin{theorem}\label{maintheorem}
Let $G$ be a finite group of exponent $n$ and let $F$ be a number field containing
a primitive $n$-th root of unity $\zeta_{n}$. Let $S$ be any finite set of finite primes of $F$.
Then there exists a finite extension $L/F$ unramified at the primes of $S$ such that for any 
$\mathcal{O}_{F}$-order $\mathcal{A}$ with $\mathcal{O}_{F}[G] \subseteq \mathcal{A} \subseteq F[G]$, 
the natural map between locally free class groups 
\[ 
\Cl(\mathcal{A}) \rightarrow \Cl(\mathcal{A} \otimes_{\mathcal{O}_{F}} \mathcal{O}_{L}), \quad
[M] \mapsto [M \otimes_{\mathcal{O}_{F}} \mathcal{O}_{L}] 
\]
is trivial.
\end{theorem}

\begin{corollary}\label{arith-disjoint}
Assume the notation and setting of Theorem \ref{maintheorem}. Let $K$ be any Galois extension of $F$ with $\Gal(K/F)=G$ such that 
\begin{enumerate}
\item $K/F$ is unramified at all finite primes outside $S$;
\item $K/F$ is linearly disjoint from $L/F$; and
\item $\mathcal{O}_{K}$ is locally free over the associated order $\mathcal{A}_{K/F}$.
\end{enumerate}
Then $\mathcal{O}_{LK}=\mathcal{O}_{K} \otimes_{\mathcal{O}_{F}} \mathcal{O}_{L}$ is free over 
$\mathcal{A}_{LK/L} = \mathcal{A}_{K/F} \otimes_{\mathcal{O}_{F}} \mathcal{O}_{L}$.
\end{corollary}

\begin{remark}
It must be noted that a capitulation result for finite-dimensional semisimple $F$-algebras can be deduced from work of Jacobinski (see \cite[Satz 7]{jacobinski}); the authors are indebted to the referee of an earlier less general version of this paper for bringing this to their attention. 
However, the crucial difference with Theorem \ref{maintheorem} is 
that Jacobinski's result does not give any control over ramification.
Also, the construction of the capitulating extension given here requires fewer steps 
than that given in \cite{jacobinski}.
\end{remark}

\begin{remark}
There are infinitely many possible choices for $L$ in 
Theorem \ref{maintheorem}. In fact, by enlarging the set $S$ appropriately, one can 
construct an infinite collection of possible choices for $L$ in which 
(the normal closures of) any two elements only have intersection equal to some extension 
of $F$ that is unramified at all finite primes.
\end{remark}

\begin{remark}
Condition (ii) of Corollary \ref{arith-disjoint} follows from condition (i) if we further assume
that for every $K'$ with $F \subsetneq K' \subseteq K$, the extension
$K'/F$ is ramified at one or more finite primes. In particular, this is true for all
extensions $K/F$ with Galois group $G$, provided $G$ is soluble, $F$ is totally complex 
(which must be the case if $n > 2$)
and the class number of $F$ is relatively prime to $n$.
\end{remark}

\begin{remark}
Let $K/F$ be a finite Galois extension of number fields with Galois group of exponent $n$.
If $F$ does not contain a primitive $n$-th root of unity $\zeta_{n}$, 
it is possible to first pass to the extension $K(\zeta_{n})/F(\zeta_{n})$ and 
then (assuming that the appropriate hypotheses hold) apply Theorem \ref{maintheorem} 
and Corollary \ref{arith-disjoint} so that there exists an extension $L/F(\zeta_{n})$ such that 
$\mathcal{O}_{LK}$ is free over $\mathcal{A}_{LK/L}$. However, by adjoining $\zeta_{n}$
it may be impossible to avoid a situation in which $L/F$ and $K/F$ are not arithmetically 
disjoint or even linearly disjoint, and so one could have 
$[LK:L]=[K(\zeta_{n}):F(\zeta_{n})] < [K:F]$.
\end{remark}

\begin{remark}
There is no ``arithmetically disjoint capitulation'' for finite Galois extensions of $p$-adic fields.
This is stated for the abelian case in \cite[Proposition 1(b)]{lettl-local}, but the proof given is also
valid for the non-abelian case. A straightforward argument then shows the following global 
statement: if $K/F$ is a finite Galois extension of number 
fields such that $\mathcal{O}_{K}$ is \emph{not} locally free over the associated order 
$\mathcal{A}_{K/F}$, then there exists no extension $L/F$ arithmetically disjoint from $K/F$ 
such that 
$\mathcal{O}_{LK}=\mathcal{O}_{K} \otimes_{\mathcal{O}_{F}} \mathcal{O}_{L}$ is free over 
$\mathcal{A}_{LK/L} = \mathcal{A}_{K/F} \otimes_{\mathcal{O}_{F}} \mathcal{O}_{L}$.
\end{remark}

\begin{remark}
By Noether's Theorem (see \cite{noether}), $\mathcal{O}_{K}$ is locally free over $\mathcal{A}_{K/F}$ 
at all primes at which $K/F$ is at most tamely ramified. Hence as a special case of
Corollary \ref{arith-disjoint}, we have a sharper version of Ichimura's capitulation 
result for relative normal integral bases of abelian extensions (see \cite{ichimura}). For primes at which $K/F$ is wildly ramified,
it is not so straightforward to determine whether we have freeness for the corresponding extension of 
$p$-adic fields. However, Lettl has shown in \cite{lettl-local} that we always have freeness for absolutely abelian extensions of $p$-adic fields. For more on the general $p$-adic case, see \cite{childs}. 
\end{remark}

\begin{remark}
There is an analogous version of Corollary \ref{arith-disjoint} for the Galois module structure of
$G$-invariant ideals of $\mathcal{O}_{L}$. For such an ideal $I$, we consider its structure as a module over
$ \mathcal{A}_{K/F}(I) := \{ x \in F[G] \, \mid \, x(I) \subseteq I \}$
and apply Theorem \ref{maintheorem} as before.
\end{remark}

\begin{remark}
With a little extra work, one can formulate and prove a similar capitulation result for
$\mathcal{O}_{F}$-orders in arbitrary semisimple $F$-algebras of finite dimension.
However, the main arithmetic application is in the special case of group algebras
where the hypotheses on $F$ in Theorem \ref{maintheorem} are easier to handle; 
in particular, it is easy to see which roots of unity must be contained in $F$.
\end{remark}

\begin{proof}[Proof of Corollary \ref{arith-disjoint}]
Conditions (i) \& (ii) and the fact that $L/F$ is unramified at the primes of $S$ imply 
that $\mathcal{O}_{L}$ and $\mathcal{O}_{K}$ are arithmetically disjoint over 
$\mathcal{O}_{F}$, i.e. 
$\mathcal{O}_{LK}=\mathcal{O}_{K} \otimes_{\mathcal{O}_{F}} \mathcal{O}_{L}$ 
(see \cite[III.2.13]{ft}). From this it is straightforward to show that 
$\mathcal{A}_{LK/L} = \mathcal{A}_{K/F} \otimes_{\mathcal{O}_{F}} \mathcal{O}_{L}$. 

By Theorem \ref{maintheorem}, we see that the class of $\mathcal{O}_{LK}$ 
is trivial in $\Cl(\mathcal{A}_{LK/L})$. This implies that the locally free $\mathcal{A}_{LK/L}$-module
$\mathcal{O}_{LK}$ is stably free. Since $L$ contains a primitive $n$-th root of unity, the Wedderburn
decomposition of the group algebra $L[G]$ is of the form $\oplus_{i=1}^{s} \Mat_{r_{i}}(L)$ and so $L[G]$
satisfies the Eichler condition relative to $\mathcal{O}_{L}$ (see \cite[Definitions 34.3 and 38.1]{reiner}). 
Therefore a result of Jacobinski (see \cite[Theorem 38.2]{reiner}, for example) shows that
$\mathcal{O}_{LK}$ is in fact free over $\mathcal{A}_{LK/L}$.
\end{proof}

\section{Capitulation for Maximal Orders}

We first prove a special case of Theorem \ref{maintheorem} which will later allow us to perform 
an important reduction step in its proof. We require the following result.

\begin{proposition}\label{maxform}
Let $r \in \mathbb{N}$, let $F$ be a number field and
let $\mathcal{O}=\mathcal{O}_{F}$ denote the ring of integers of $F$.
For each ideal $\mathfrak{a}$ of $\mathcal{O}$, let
$$ \Lambda_{\mathfrak{a},r} =
\left(\begin{array}{cccc}
\mathcal{O} & \cdots & \mathcal{O}  & \mathfrak{a}^{-1} \\
\vdots & \ddots & \vdots & \vdots  \\
\mathcal{O} & \cdots & \mathcal{O} & \mathfrak{a}^{-1}  \\
\mathfrak{a} & \cdots & \mathfrak{a}  & \mathcal{O}
\end{array}\right)$$
denote the ring of all $r \times r$ matrices $(x_{ij})$ where $x_{11}$ ranges over all elements of 
$\mathcal{O}$, \ldots, $x_{1r}$ ranges over all elements of $\mathfrak{a}^{-1}$, and so on. 
(For $r=1$, we take $\Lambda_{\mathfrak{a},r}=\mathcal{O}$.)
Then $\Lambda_{\mathfrak{a},r}$ is a maximal $\mathcal{O}$-order in $\Mat_{r}(F)$ and every maximal 
$\mathcal{O}$-order in $\Mat_{r}(F)$ is isomorphic to one
of this form, for some ideal $\mathfrak{a}$ of $\mathcal{O}$.
\end{proposition}

\begin{proof}
This is a special case of \cite[Corollary 27.6]{reiner}.
\end{proof}

\begin{proposition}\label{maxcap}
Let $G$ be a finite group of exponent $n$ and let $F$ be a number field containing
a primitive $n$-th root of unity $\zeta_{n}$.
Let $H$ be the Hilbert class field of $F$ and let $\mathcal{M}_{F[G]}$ be any maximal 
$\mathcal{O}_{F}$-order of $F[G]$.
Then $\mathcal{M}_{H[G]} := \mathcal{M}_{F[G]} \otimes_{\mathcal{O}_{F}} \mathcal{O}_{H}
\cong \oplus_{i=1}^{s} \Mat_{r_{i}}(\mathcal{O}_{H})$ for some $r_{i},s \in \mathbb{N}$,
and so is a maximal $\mathcal{O}_{H}$-order of $H[G]$. 
Furthermore, the natural map between locally free class groups 
\[ 
\theta : \Cl(\mathcal{M}_{F[G]}) \rightarrow \Cl(\mathcal{M}_{H[G]}), 
\quad [M] \mapsto [M \otimes_{\mathcal{O}_{F}} \mathcal{O}_{H}] 
\]
is trivial.
\end{proposition}

\begin{proof}
As $F$ and $H$ both contain a primitive $n$-th root of unity, 
the Wedderburn decompositions of the group algebras
$F[G]$ and $H[G]$ are of the forms $\oplus_{i=1}^{s} \Mat_{r_{i}}(F)$ and $\oplus_{i=1}^{s} \Mat_{r_{i}}(H)$,
respectively. Henceforth abbreviate $\mathcal{M}_{F[G]}$ and $\mathcal{M}_{H[G]}$ to $\mathcal{M}$
and $\mathcal{M}'$, respectively.
By Proposition \ref{maxform}, we have
$\mathcal{M} \cong \oplus_{i=1}^{s} \Lambda_{\mathfrak{a}_{i},r_{i}}$ 
for some ideals $\mathfrak{a}_{i}$ of $\mathcal{O}_{F}$. Identifying $\mathcal{M}$ with  
$\oplus_{i=1}^{s} \Lambda_{\mathfrak{a}_{i},r_{i}}$, let $\mathcal{M}_{i} = \Lambda_{\mathfrak{a}_{i},r_{i}}$
and $\mathcal{M}_{i}' = \mathcal{M}_{i} \otimes_{\mathcal{O}_{F}} \mathcal{O}_{H} = \Lambda_{(\mathfrak{a}_{i}\mathcal{O}_{H}),r_{i}}$. Then $\mathcal{M}' = \oplus_{i=1}^{s} \mathcal{M}_{i}'$ and since 
each ideal $\mathfrak{a}_{i}\mathcal{O}_{H}$ is principal, we have 
$\mathcal{M}_{i}' \cong \Mat_{r_{i}}(\mathcal{O}_{H})$. Hence 
$\mathcal{M}_{i}'$ is a maximal $\mathcal{O}_{H}$-order of $\Mat_{r_{i}}(H)$ for each $i$
by Proposition \ref{maxform},
and so $\mathcal{M}'$ is a maximal $\mathcal{O}_{H}$-order of $H[G]$.

Showing that $\theta$ is trivial is equivalent to showing that each component 
\[ 
\theta_{i} : \Cl(\mathcal{M}_{i}) \rightarrow \Cl(\mathcal{M}_{i}'), 
\quad [N] \mapsto [N \otimes_{\mathcal{O}_{F}} \mathcal{O}_{H}] 
\]
is trivial. Fix $i$ and let $N_{i}$ be a locally free $\mathcal{M}_{i}$-module (of rank $t$ for 
some $t \in \mathbb{N}$). 
It suffices to show that $N_{i} \otimes_{\mathcal{O}_{F}} \mathcal{O}_{H}$ is a free 
$\mathcal{M}_{i}'$-module (of rank $t$).

We shall now assume some basic facts on Morita equivalence (\cite[Chapter 4]{reiner} contains the relevant material; note for later use that for every commutative ring $S$ and $n\in \mathbb{N}$, the matrix ring $\Mat_n(S)$ is Morita equivalent to $S$, since it is the $S$-endomorphism ring of the free module $S^n$). For a ring $R$, let $\Mod(R)$ denote the category of finitely generated projective $R$-modules.
We have the following diagram
\[
\xymatrix@1@!0@=48pt { 
\Mod(\mathcal{O}_{F})  \ar@{>}[d]_{\beta} \ar@{<-}[rr]^{\alpha} & 
& \Mod(\mathcal{M}_{i})  \ar@{>}[d]^{\delta} \\
\Mod(\mathcal{O}_{H}) \ar@{->}[rr]^{\gamma} &
& \Mod(\mathcal{M}_{i}')
}
\]
(where $\alpha$ and $\gamma$ are the Morita functors, and $\beta$ and $\delta$ are both the
$- \otimes_{\mathcal{O}_{F}} \mathcal{O}_{H}$ functor), which commutes up to natural equivalence. 
First note that $N_{i}$ is $\mathcal{M}_{i}$-projective, locally free of rank $t$.
Thus $\alpha(N_{i})$ is $\mathcal{O}_{F}$-projective, locally free of rank $r_{i}t$ and 
so $\beta(\alpha(N_{i}))$ is $\mathcal{O}_{H}$-free of rank $r_{i}t$. Finally, since 
$\beta(\alpha(N_{i})) \cong \oplus_{j=1}^{r_{i}t} \mathcal{O}_{H}$ 
and $\mathcal{M}_{i}' \cong \Mat_{r_{i}}(\mathcal{O}_{H})$, the Morita functor $\gamma$ maps
$\beta(\alpha(N_{i}))$ to a module isomorphic to $\oplus_{j=1}^{t} \Mat_{r_{i}}(\mathcal{O}_{H})$, and so
$\delta(N_{i}) \cong \gamma(\beta(\alpha(N_{i})))$ is a free $\mathcal{M}_{i}'$-module of rank $t$.
\end{proof}

\section{Capitulation for Units}\label{cap-for-units}

In this section, we prove a ``capitulation result for units''. The key idea in the proof of Theorem 
\ref{maintheorem} is to combine this with the fact that the kernel group of the locally free class group 
can be described in terms of units.

The following uses and develops some ideas from
\cite[Lemma 8]{ichimura} and \cite[Lemma 3.1]{repunits}.

\begin{proposition}
Let $F$ be a number field and let $m$ be any positive integer. Let $u \in \mathcal{O}_{F}$ be any element with $(u,m)=1$. Then there exists a finite extension $L/F$ unramified at all prime divisors of $m$ such that there exists $\varepsilon \in \mathcal{O}_{L}^{\times}$ with $\varepsilon \equiv u \, \mod \, m$.
\end{proposition}

\begin{proof}
The hypothesis that $(u,m)=1$ shall be used throughout without any further mention.

Let $t \in \mathbb{N}$ such that $u^{t} \equiv 1 \, \mod \, m^{2}$.
Without loss of generality, we can assume that $m$ divides $t$ and $t \geq m+4$.
There exist elements $a,b \in \mathcal{O}_{F}$ such that
\[
au^{t-1} + bm^{t} = \frac{1-u^{t}}{m}\, .
\]
Putting $c=a-m^{t-1}$ and $d=bm+u^{t-1}$, we have
\begin{equation}\label{relprimeeqn}
cu^{t-1} + dm^{t-1} = \frac{1-u^{t}}{m}\, .
\end{equation}
We see that $m$ divides $c$ because $m$ divides both 
$dm^{t-1}$ and $(1-u^{t})/m$. Furthermore,
$(d,m)=1$ since $d=bm+u^{t-1}$. 

Let $\theta=\theta_{1}, \theta_{2}, \ldots, \theta_{t}$ be the roots of the polynomial
\[
f(X) = X^{t} + cX^{t-1} -mX^{m+2} + uX^{m+1} + d \in \mathcal{O}_{F}[X]
\]
and let $L=F(\theta)$. 
Then we have
\begin{eqnarray*}
\prod_{i=1}^{t} (u - m\theta_{i}) &=& 
m^{t} f(u/m) = u^{t} + mcu^{t-1} - mm^{t-(m+2)}u^{m+2} + um^{t-(m+1)}u^{m+1} + dm^{t} \\
&=&  u^{t} + mcu^{t-1} + dm^{t} =  u^{t} + m(cu^{t-1} + dm^{t-1}) = u^{t} + (1-u^{t}) =1,
\end{eqnarray*}
where the penultimate equality is due to (\ref{relprimeeqn}).
Therefore, $\varepsilon = u - m \theta \in \mathcal{O}_{L}^{\times}$ and 
$\varepsilon \equiv u \, \mod \, m$.

Now let $\mathfrak{p}$ be a prime of $\mathcal{O}_F$ dividing $m$.
The formal derivative of $f(X)$ is 
\[
f'(X) = tX^{t-1} + c(t-1)X^{t-2} - (m+2)mX^{m+1} + (m+1)uX^{m}.
\]
However, $\mathfrak{p}$ divides $m$, and $m$ divides both $t$ and $c$, so we have
\[
f'(X) \equiv uX^{m} \, \mod \, \mathfrak{p} \, .
\]
Note that $u \not \equiv 0 \, \mod \, \mathfrak{p}$. Hence 
$X=0$ is the unique root of $f'(X) \, \mod \, \mathfrak{p}$. However,
$f(0) = d \not \equiv 0 \, \mod \, \mathfrak{p}$ since $(d,m)=1$. Therefore $f(X)$ and $f'(X)$ have no common roots $\mod \, \mathfrak{p}$ and so $L/F$ is unramified at $\mathfrak{p}$.
\end{proof}

\begin{remark}
The construction of $L$ in the above proof is far from canonical. An interesting but apparently
difficult question is whether $L/F$ can always be taken to be relatively abelian.
\end{remark}

\begin{corollary}\label{null-units}
Let $F$ be a number field and let $m$ be any positive integer. Then there exists a finite extension $L/F$ unramified at all prime divisors of $m$ such that the natural map
\[
\frac{(\mathcal{O}_{F}/m\mathcal{O}_{F})^{\times}}{\phi(\mathcal{O}_{F}^{\times})}
\longrightarrow \frac{(\mathcal{O}_{L}/m\mathcal{O}_{L})^{\times}}{\psi(\mathcal{O}_{L}^{\times})}\, ,
\]
with $\phi:\mathcal{O}_{F}^{\times} \rightarrow (\mathcal{O}_{F}/m\mathcal{O}_{F})^{\times}$
and  $\psi:\mathcal{O}_{L}^{\times} \rightarrow (\mathcal{O}_{L}/m\mathcal{O}_{L})^{\times}$
the natural projections, is trivial.
\end{corollary}

\section{Proof of the Main Theorem}

We now bring together the results of the previous two sections to prove the main theorem.

\begin{proof}[Proof of Theorem \ref{maintheorem}]
Let $\mathcal{M}_{F[G]}$ be a maximal $\mathcal{O}_{F}$-order containing $\mathcal{A}$
and let $H$ be the Hilbert class field of $F$.
Let $\mathcal{B}= \mathcal{A} \otimes_{\mathcal{O}_{F}} \mathcal{O}_{H}$ and 
$\mathcal{M}_{H[G]} = \mathcal{M}_{F[G]} \otimes_{\mathcal{O}_{F}} \mathcal{O}_{H}$.
Note that $\mathcal{M}_{H[G]}$ is a maximal $\mathcal{O}_{H}$-order in $H[G]$ by 
Proposition \ref{maxcap}.
Since $\mathcal{O}_{H}$ is projective over
$\mathcal{O}_{F}$, the functor $- \otimes_{\mathcal{O}_{F}} \mathcal{O}_{H}$ is exact. Hence
we have the following commutative diagram
\[
\xymatrix@1@!0@=48pt { 
0 \ar@{->}[rr] & & \D(\mathcal{A}) \ar@{->}[rr] \ar@{->}[d]_{\alpha} & & \Cl(\mathcal{A}) \ar@{->}[rr]^{\delta}
\ar@{->}[d]_{\beta} & & \Cl(\mathcal{M}_{F[G]})  \ar@{->}[d]_{\gamma} \\
0 \ar@{->}[rr] & & \D(\mathcal{B}) \ar@{->}[rr] 
& & \Cl(\mathcal{B}) \ar@{->}[rr]^{\varepsilon} 
& & \Cl(\mathcal{M}_{H[G]}) 
}
\]
where $\D(\mathcal{A}):= \ker \delta$, $\D(\mathcal{B}):= \ker \varepsilon$,
and $\alpha, \beta$ and $\gamma$ are the maps induced by $- \otimes_{\mathcal{O}_{F}} \mathcal{O}_{H}$. 

The map $\gamma$ is trivial by Proposition \ref{maxcap}, and so we have 
$\beta(\Cl(\mathcal{A})) \subseteq  \D(\mathcal{B})$. Hence we are reduced to showing that there exists 
an extension $L/H$ unramified at all primes in $S$ such that natural map
\[
\D(\mathcal{B}) \longrightarrow \D(\mathcal{B} \otimes_{\mathcal{O}_{H}} \mathcal{O}_{L}) 
\subseteq \Cl(\mathcal{B} \otimes_{\mathcal{O}_{H}} \mathcal{O}_{L}), 
\quad [M] \mapsto [M \otimes_{\mathcal{O}_{H}} \mathcal{O}_{L}]
\]
is trivial.

We henceforth abbreviate $\mathcal{M}_{H[G]}$ to $\mathcal{M}$. 
By Proposition \ref{maxcap}, we have 
$\mathcal{M} \cong \bigoplus_{i=1}^{s} \Mat_{r_{i}}(\mathcal{O}_{H})$ for some $r_{i},s \in \mathbb{N}$.  
Let $m \in \mathbb{N}$ be a multiple of $|G|$ that is divisible by all the primes in $S$. Note that $m\mathcal{M} \subseteq \mathcal{B}$.

We now assume some basic facts from K-theory (for an introduction, see \cite{curtis-reiner2}). 
We have the following Milnor square
\[
\xymatrix@1@!0@=48pt { 
\mathcal{B} \ar@{>}[d] \ar@{->}[rr] & & \mathcal{M}  \ar@{>}[d] \\
\mathcal{B}/m\mathcal{M}\ar@{->}[rr] & & \mathcal{M}/m\mathcal{M}
}
\]
where the horizontal arrows are the natural inclusions and the vertical arrows are the natural projections
(note that this is a special case of a fiber product).
By \cite[p. 242]{curtis-reiner2} (with 
$\Lambda = \mathcal{B}, \Gamma= \mathcal{M}, \overline{\Lambda}=\mathcal{B}/m\mathcal{M}$ 
and $\overline{\Gamma}=\mathcal{M}/m\mathcal{M}$) 
we have the following exact sequence
\[
\K_{1}(\mathcal{M}) \times \K_{1}(\mathcal{B}/m\mathcal{M})
\stackrel{\phi}\longrightarrow 
\K_{1}(\mathcal{M}/m\mathcal{M})
\stackrel{\partial}\longrightarrow
\D(\mathcal{B})
\longrightarrow 0
\]
where $\partial$ is given by the Milnor patching process and is described explicitly in the 
proof of Milnor's Theorem \cite[Theorem 42.13]{curtis-reiner2}. By Morita equivalence, this becomes
\[
\K_{1}(\oplus_{i=1}^{s} \mathcal{O}_{H}) \times \K_{1}(\mathcal{B}/m\mathcal{M}) 
\stackrel{\phi}\longrightarrow 
\K_{1}(\oplus_{i=1}^{s} (\mathcal{O}_{H} / m \mathcal{O}_{H})) 
\stackrel{\partial}\longrightarrow 
\D(\mathcal{B}) 
\longrightarrow 0,
\]
which can be rewritten as
\[
\textstyle
(\prod_{i=1}^{s} \mathcal{O}_{H}^{\times}) \times \K_{1}(\mathcal{B}/m\mathcal{M}) 
\stackrel{\phi}{\longrightarrow} 
(\prod_{i=1}^{s} (\mathcal{O}_{H} / m \mathcal{O}_{H})^{\times}) 
\stackrel{\partial}{\longrightarrow} \D(\mathcal{B}) \longrightarrow 0.
\]
Note that $\phi$ restricted to $(\prod_{i=1}^{s} \mathcal{O}_{H}^{\times})$ is just the natural projection 
and $\phi(\prod_{i=1}^{s} \mathcal{O}_{H}^{\times}) \subseteq \ker\partial$, giving a surjection
\[
g: \prod_{i=1}^{s} \frac{(\mathcal{O}_{H}/m\mathcal{O}_{H})^{\times}}{\phi(\mathcal{O}_{H}^{\times})}
\longrightarrow \D(\mathcal{B})\, .
\]
By Corollary \ref{null-units}, there exists an extension 
$L/H$ unramified at all prime divisors of $m$ (in particular, at all primes in $S$) 
such that if 
$f: \prod_{i=1}^{s} (\mathcal{O}_{H}/m\mathcal{O}_{H})^{\times}
\rightarrow \prod_{i=1}^{s} (\mathcal{O}_{L}/m\mathcal{O}_{L})^{\times}$
is the natural map, then the induced map
\[ \bar{f}:
\prod_{i=1}^{s} \frac{(\mathcal{O}_{H}/m\mathcal{O}_{H})^{\times}}{\phi(\mathcal{O}_{H}^{\times})} \longrightarrow 
\prod_{i=1}^{s} \frac{(\mathcal{O}_{L}/m \mathcal{O}_{L})^{\times}}{\psi(\mathcal{O}_{L}^{\times})}
\]
is trivial ($\phi$ and $\psi$ are the natural projections). 

The functor $ - \otimes_{\mathcal{O}_{H}} \mathcal{O}_{L}$ is exact because
$\mathcal{O}_{L}$ is projective over $\mathcal{O}_{H}$. Defining
$\mathcal{B'}=\mathcal{B} \otimes_{\mathcal{O}_{H}} \mathcal{O}_{L}$ and
$\mathcal{M'}=\mathcal{M} \otimes_{\mathcal{O}_{H}} \mathcal{O}_{L}$, 
we therefore have a second Milnor square
\[
\xymatrix@1@!0@=48pt { 
\mathcal{B'}  \ar@{>}[d] \ar@{->}[rr] & 
& \mathcal{M'}  \ar@{>}[d] \\
\mathcal{B'} /m\mathcal{M'} \ar@{->}[rr] &
& \mathcal{M'}/m\mathcal{M'}
}
\]
where again the horizontal arrows are the natural inclusions and the vertical arrows are the natural projections.
The two Milnor squares can be ``glued together'' using the base change maps 
$\mathcal{B} \rightarrow \mathcal{B'}$, etc. to form the following commutative cube.
\[
\xymatrix{ 
& \mathcal{B} \ar@{->}[rr]\ar@{->}'[d][dd] 
& & \mathcal{M} \ar@{->}[dd] 
\\ 
\mathcal{B'} \ar@{<-}[ur]\ar@{->}[rr]\ar@{->}[dd] 
& &  \mathcal{M}' \ar@{<-}[ur]\ar@{->}[dd] 
\\
& \mathcal{B} /m\mathcal{M} \ar@{->}'[r][rr] 
& &   \mathcal{M} /m\mathcal{M}
\\
\mathcal{B'} /m\mathcal{M'} \ar@{->}[rr]\ar@{<-}[ur] 
& & \mathcal{M'} /m\mathcal{M'} \ar@{<-}[ur]
}
\]
Moreover, the Milnor patching process can easily be shown to commute with base change.
Hence we have a commutative square
\[
\xymatrix@1@!0@=48pt { 
\prod_{i=1}^{s}(\mathcal{O}_{H}/m\mathcal{O}_{H})^{\times}
 \ar@{>}[d]_{f} \ar@{->}[rr]^{\qquad \partial} & &
\D(\mathcal{B}) \ar@{>}[d]^{h} \\
\prod_{i=1}^{s}(\mathcal{O}_{L}/m\mathcal{O}_{L})^{\times}
\ar@{->}[rr]^{\qquad \partial'} & & 
\D(\mathcal{B'})}
\]
where $\partial'$ is defined analogously to $\partial$ and
$h$ is the map $[M] \mapsto [M \otimes_{\mathcal{O}_{H}} \mathcal{O}_{L}]$.
Consequently, the following diagram commutes as well.
\[
\xymatrix@1@!0@=48pt { 
\prod_{i=1}^{s} \frac{(\mathcal{O}_{H}/m\mathcal{O}_{H})^{\times}}{\phi(\mathcal{O}_{H}^{\times})}
 \ar@{>}[d]_{\bar{f}} \ar@{->}[rr]^{\qquad g} & &
\D(\mathcal{B}) \ar@{>}[d]^{h} \\
\prod_{i=1}^{s} \frac{(\mathcal{O}_{L}/m\mathcal{O}_{L})^{\times}}{\psi(\mathcal{O}_{L}^{\times})}
\ar@{->}[rr]^{\qquad g'} & & 
\D(\mathcal{B'})}
\]
As $g$ is surjective and $\bar{f}$ is trivial, this implies that $h$ is trivial, 
which is exactly what we wanted to show.
\end{proof}

\section{Non-Capitulation}

As a complement to our results on capitulation, we make the following easy observation.

\begin{proposition}\label{non-cap}
Let $F$ be a number field, $G$ be a finite group and
$\mathcal{A}$ be any $\mathcal{O}_{F}$-order with 
$\mathcal{O}_{F}[G] \subseteq \mathcal{A} \subseteq F[G]$.
Let $N$ be a finitely generated locally free $\mathcal{A}$-module 
and denote the order of the class $[N]$ in $\Cl(\mathcal{A})$ by $o([N])$.
Suppose that $L/F$ is a finite extension
such that $o([N]) \nmid [L:F]$ and $\mathcal{O}_{L}$ is free over $\mathcal{O}_{F}$. 
Then $[N]$ is not in the kernel of the map
\[
\Cl(\mathcal{A}) \rightarrow  \Cl(\mathcal{A} \otimes_{\mathcal{O}_{F}} \mathcal{O}_{L}), \quad
[M] \mapsto [M \otimes_{\mathcal{O}_{F}} \mathcal{O}_{L}] \, .
\]
\end{proposition}

\begin{proof}
Suppose that $[N \otimes_{\mathcal{O}_{F}} \mathcal{O}_{L}]$ is trivial in $\Cl(\mathcal{A} \otimes_{\mathcal{O}_{F}} \mathcal{O}_{L})$. Then $N \otimes_{\mathcal{O}_{F}} \mathcal{O}_{L}$ is stably free over 
$\mathcal{A} \otimes_{\mathcal{O}_{F}} \mathcal{O}_{L}$ and so there exist $n,m \in \mathbb{N}$ such that
\[
(\mathcal{A} \otimes_{\mathcal{O}_{F}} \mathcal{O}_{L})^{n} \oplus (N \otimes_{\mathcal{O}_{F}} \mathcal{O}_{L}) \cong 
(\mathcal{A} \otimes_{\mathcal{O}_{F}} \mathcal{O}_{L})^{n+m}
\textrm{ as } (\mathcal{A} \otimes_{\mathcal{O}_{F}} \mathcal{O}_{L}) \textrm{-modules.}
\]
Then restricting coefficients gives
\[
(\mathcal{A} \otimes_{\mathcal{O}_{F}} \mathcal{O}_{L})^{n} \oplus (N \otimes_{\mathcal{O}_{F}} \mathcal{O}_{L}) \cong 
(\mathcal{A} \otimes_{\mathcal{O}_{F}} \mathcal{O}_{L})^{n+m}
\textrm{ as } \mathcal{A} \textrm{-modules.} \]
Since $\mathcal{O}_{L}$ is free over $\mathcal{O}_{F}$, we have
\[
\mathcal{A}^{n[L:F]} \oplus N^{[L:F]} \cong 
\mathcal{A}^{(n+m)[L:F]}
\textrm{ as } \mathcal{A} \textrm{-modules.}
\]
Hence $[N]^{[L:F]}$ is trivial in $\Cl(\mathcal{A})$, contradicting the hypothesis that $o([N]) \nmid [L:F]$.
\end{proof}

\begin{corollary}
Let $F$, $G$ and $\mathcal{A}$ be as in Proposition \ref{non-cap}.
Let $L/F$ be a finite extension such that $([L:F],|\Cl(\mathcal{A})|)=1$ 
and $\mathcal{O}_{L}$ is free over $\mathcal{O}_{F}$.
Then the map
\[
\Cl(\mathcal{A}) \rightarrow  \Cl(\mathcal{A} \otimes_{\mathcal{O}_{F}} \mathcal{O}_{L}), \quad
[M] \mapsto [M \otimes_{\mathcal{O}_{F}} \mathcal{O}_{L}]
\]
is injective.
\end{corollary}

\section{Acknowledgements}

The authors are grateful to the Deutscher Akademischer Austausch Dienst (German Academic Exchange Service) for a grant allowing the second named author to visit the first for the 2006-07 academic year, thus making this collaboration possible. The authors are indebted to both the referee
of an earlier less general version of this paper and the referee of the present version
for numerous helpful comments and suggestions.

\end{document}